\RequirePackage{fix-cm}
\documentclass[a4paper, 12pt, leqno]{article}%{amsart}

\usepackage{latexsym, amsmath, amsthm, amssymb, graphicx, graphics, epsfig, bm}
\usepackage{bm}
\usepackage{tikz-cd}

\usepackage{enumitem}

\swapnumbers 
\usepackage{xpatch}
\xpatchcmd\swappedhead{~}{.~}{}{}

\tikzset{cong/.style={draw=none,edge node={node [sloped, allow upside down, auto=false]{$\cong$}}},
         Isom/.style={every to/.append style={edge node={node [above,sloped, inner sep=0.4pt, allow upside down, auto=false]{$\sim$}}}}}

\usepackage{fullpage}

\newcommand{\comment}[1]{}

\usepackage[T1]{fontenc}
\usepackage{latexsym}
\usepackage{manfnt}
\usepackage[sans]{dsfont}
\usepackage{palatino}
\usepackage[sc, osf]{mathpazo} 
\usepackage{eulervm}
\usepackage{euscript}
\newcommand{\sub}[1]{\vspace{8pt}\textbf{#1}\hspace*{0.5em}}
\newtheorem{thm}{Theorem}%[section]
\newtheorem*{thm*}{Theorem}
\newtheorem{lem}[thm]{Lemma}
\newtheorem*{lem*}{Lemma}

\newtheorem*{cor*}{Corollary}

\newtheorem*{prop*}{Proposition}
\newtheorem*{claim*}{Claim}

\theoremstyle{definition}

\renewcommand{\proof}{\vspace{-8pt}\noindent\textit{\textbf{Proof. }}}

\renewcommand{\endproof}{$\square$}

\renewcommand{\leq}{\leqslant}

\newcommand{\set}[1]{\left\{#1\right\}}

\newcommand{\norm}[1]{\left\|#1\right\|}

\newcommand{\inv}{^{-1}}

\newcommand{\eps}{\varepsilon}
\renewcommand{\phi}{\varphi}

\newcommand{\R}{\mathbb R}

\renewcommand{\H}{\mathbb H}
\newcommand{\A}{\mathbb A}

\DeclareMathOperator{\im}{Im}

\newcommand{\xX}{x\in X}

\newcommand{\f}{\gamma} %{\bm{\chi}}

\title{The continuity properties of discrete-spectrum families of Fredholm operators} 
\author{Marina\,Prokhorova}
\date{}

\begin{document}

\maketitle

\footnotetext{Department of Mathematics, Technion -- Israel Institute of Technology}
\footnotetext{This work was partially supported by ISF grants no. 431/20 and 844/19}

%\vspace*{-1ex}

\renewcommand{\baselinestretch}{1.02}
\selectfont

\sub{Introduction.}
Recently, N. Ivanov  \cite{i1} introduced notions of \emph{Fredholm families} and 
\emph{discrete-spectrum} families of self-adjoint Fredholm operators.
Fredholm families are not (and probably cannot be) defined in terms of continuity in any topology. 
The discrete-spectrum families are introduced in \cite{i1} as a natural analogue of Fredholm
families for self-adjoint operators with discrete spectrum
and are defined in \cite{i1} in terms of Fredholm families.

The goal of this note is to relate the notion of discrete-spectrum families with the classical continuity properties.
We show that the discrete-spectrum families are exactly the families of
self-adjoint operators with discrete spectrum which are continuous in the uniform resolvent topology
and clarify Ivanov's approach \cite{i1, i2} to a theorem of Melrose--Piazza \cite{mp}.
See Theorems \ref{thm:1} and \ref{thm:2} below, which are the main results of this note.

I am grateful to N. Ivanov for asking me to write this paper.

\sub{Fredholm and discrete-spectrum families.}
Let us recall some definitions from \cite{i1, i2}.
All Hilbert spaces are assumed to be complex separable and infinite-dimensional.
Let $\H$ be a locally trivial Hilbert bundle over $X$ with the fibers $H_x$, $x \in X$.
Let $\A$ be a family of self-adjoint operators $A_x\colon H_x\to H_x$. 
We assume either that operators $A_x$ are bounded or that they are unbounded and closed densely defined.
We also assume, following \cite{as} and \cite{i1, i2}, 
that operators $A_x$ are neither essentially positive nor essentially negative.
%The assumption about the essential spectrum  
The last assumption is important for the index theory, but not for our arguments.

A pair $(U,\eps)$, where $U\subset X$ is an open set and $\eps>0$, 
is said to be \emph{adapted} to the family $\A$ 
if for every $x\in U$ 
the image $V_x = \im P_{[-\eps,\eps]}(A_x)$ of the spectral projection is finite-dimensional, 
$\pm\eps$ does not belong to the spectrum $\sigma(A_x)$ of $A_x$, %$\pm \eps\notin\sigma(A_x)$,  
and both $V_x$ and the restriction of $A_x$ to $V_x$ depend norm continuously on $x\in U$.

The family $\A$ is said to be a \emph{Fredholm family} if all operators $A_x$ are Fredholm 
and for every $x\in X$ there exists a pair $(U,\eps)$ adapted to $\A$ such that $x \in U$. 

The family $\A$ is said to be a \emph{discrete-spectrum family} if for every $\lambda\in\R$ 
the family of operators $A_x-\lambda$ is a Fredholm family.
Clearly, every operator $A_x$ in a discrete-spectrum family has discrete spectrum
and hence is an unbounded operator with compact resolvent.

\begin{lem}\label{lem:1}
A family $\A$ of self-adjoint operators is a discrete-spectrum family if and only if 
for every $x\in X$ and $b>0$ there is a pair $(U,c)$ adapted to $\A$ such that $x\in U$ and $c>b$.
\end{lem}

\proof
Suppose that $\A$ is a discrete-spectrum family.
%the condition from the lemma holds.
Let us fix an arbitrary $\lambda\in\R$ and 
choose a pair $(U,c)$ adapted to $A$ such that $c>|\lambda|$.
Then there exists some $\eps>0$ such that $\lambda\pm\eps\in[-c,c]\setminus\sigma(A_x)$.
By the standard properties of self-adjoint operators in a finite-dimensional space, 
there is a neighborhood $U'\subset U$ of $x$ such that $(U',\eps)$ is adapted to $\A-\lambda$.
This prove the ``if'' part of the lemma.

Let us prove the ``only if'' part.
Let $x\in X$ and $b>0$.
Since $\A\pm b$ are Fredholm families, there exists
$c>b$ such that $\pm c\notin\sigma(A_x)$.
By the definition of a discrete-spectrum family, 
for every $\lambda\in\R$ there exists a pair $(U_{\lambda},\eps_{\lambda})$ adapted to $\A-\lambda$.
By compactness, there exists a finite subset $\Lambda\subset\R$ such that 
the intervals $(\lambda-\eps_{\lambda},\lambda+\eps_{\lambda})$ 
cover $[-c,c]$ when $\lambda$ runs over $\Lambda$. 
Let 
\[
c^- = \min_{\lambda\in\Lambda}\set{\lambda-\eps_{\lambda}}
\quad\mbox{and}\quad 
c^+ = \max_{\lambda\in\Lambda}\set{\lambda+\eps_{\lambda}}.
\]
Then $c^-<-c<0<c<c^+$ and
the restriction $A'_y$ of $A_y$ to $\im P_{[c^-,\,c^+]}(A_y)$ 
depends norm continuously on 
\[
y\in U' = \bigcap\nolimits_{\lambda\in\Lambda} U_{\lambda}.
\]
Since $\pm c\notin\sigma(A'_x)$, 
by the standard properties of finite rank operators 
they are also outside of the spectrum of $A'_y$ for $y$ 
belonging to some neighborhood $U\subset U'$ of $x$.
Then $(U,c)$ is adapted to $\A$.
This proves the ``only if'' part.
\endproof

\sub{Graph (uniform resolvent) topology.}
Recall that the \emph{graph topology} on the space of self-adjoint 
(closed densely defined) operators coincides with the \emph{uniform resolvent topology},
that is, the topology induced by the map $A\mapsto(A+i)\inv$ 
from the norm topology on bounded operators.

\begin{thm}\label{thm:1}
A family $\A$ of self-adjoint operators is a discrete-spectrum family if and only if
$\A$ is a graph continuous family of operators with compact resolvent.
\end{thm}

\proof
Since the statement is local, we can fix a local trivialization of $\H$ 
and consider all operators as acting in the same Hilbert space $H$.

If $\A$ is a graph continuous family of operators with compact resolvent, 
then $\A-\lambda$ has the same property for every $\lambda\in\R$.
Since every graph continuous family of operators with compact resolvent 
is obviously Fredholm, this proves the ``if'' part.

Let us prove the ``only if'' part.
Let $\xX$ and $\delta>0$. 
By Lemma \ref{lem:1}, there is a pair $(U,c)$ adapted to $\A$ such that $c>\delta\inv$.
Let $B_y$ be the operator equal to $(A_y+i)\inv$ on $V_y = \im P_{[-c,\,c]}(A_y)$ 
and to zero on $V_y^{\bot}$.
By our choice of $c$, 
\[
\norm{(A_y+i)\inv - B_y} <\delta
\] 
for every $y\in U$.
By the definition of an adapted pair, $B_y$ continuously depends on $y\in U$.
Let $U'\subset U$ be a neighborhood of $x$ such that $\norm{B_y - B_x}<\delta$ for $y\in U'$.
Then 
\[
\norm{(A_y+i)\inv - (A_x+i)\inv} < 3\delta
\] 
for every $y\in U'$.
Since $\xX$ and $\delta>0$ were chosen arbitrarily, $\A$ is graph continuous.
This proves the ``only if'' part.
\endproof

\sub{Riesz topology.}
Let $\gamma\colon\R\to[0,1]$ be the function defined by $\gamma(t)=t(1+t^2)^{-1/2}$.
The Riesz topology on self-adjoint (closed densely defined) operators is the topology 
induced by the \emph{bounded transform} 
\[
A\longmapsto \gamma(A)= A(1+A^2)^{-1/2}
\] 
from the norm topology on the space of bounded operators.

\sub{Adapted trivializations.}
Let us recall some further definitions from \cite{i1}, \cite{i2}.
A local trivialization of $\H$ over an open subset $U\subset X$ 
is said to be \emph{strictly adapted} to the family $\A$ if 
for every $x\in U$ there exists a neighborhood $U'\subset U$ of $x$ and $\eps>0$
such that $(U',\eps)$ is adapted to $\A$ and this trivialization takes
the family of spectral projections $P_{[\eps,\infty)}(A_y)$, $y\in U'$
into a norm continuous family.
A local trivialization of $\H$ over an open subset $U\subset X$ 
is said to be \emph{fully adapted} to the family $\A$ if 
it takes $A_y$, $y\in U$ to a Riesz continuous family.
Clearly, a fully adapted trivialization is strictly adapted.

The family $\A$ is said to be a \emph{strictly Fredholm family} if it is Fredholm
and the Hilbert bundle $\H$ admits local trivializations strictly adapted to $\A$ over open subsets
forming a covering of $X$.
\emph{Fully Fredholm families} are defined similarly.

\begin{thm}\label{thm:2}
Let $\A$ be a discrete-spectrum family.
Then every local trivialization strictly adapted to $\A$ is fully adapted to $\A$.
In particular, every discrete-spectrum and strictly Fredholm family is fully Fredholm.
\end{thm}

\proof
Let us fix a strictly adapted to $\A$ trivialization of $\H$ over $U$ and use it
to identify Hilbert spaces $H_x$, $x\in U$ with $H$.
Let $x\in U$ and $\delta>0$. 
There is a pair $(U',\eps)$ adapted to $\A$ such that $x\in U'\subset U$ 
and the map $y\mapsto P_{[\eps,\infty)}(A_y)$ is norm continuous on $U'$.
Using Lemma \ref{lem:1} and decreasing $U'$ if needed, 
we get a pair $(U',c)$ adapted to $\A$ such that $x\in U'$, $c>\eps$, and $\f(c) > 1-\delta$.
By the definition of an adapted pair, 
the spectral projection $q_y = P_{(-c,\,c)}(A_y)$ and the restriction $A'_y$ of $A_y$ to the range of $q_y$ 
are norm continuous on $U'$.
It follows that the families of spectral projections 
\[
y \mapsto q^+_y = P_{[c,+\infty)}(A_y) = P_{[\eps,\infty)}(A_y) - P_{[\eps,c)}(A'_y) \quad\mbox{and} 
\]
\vspace{-16pt}
\[ 
\hspace*{-6.6em}
y \mapsto q^-_y = P_{(-\infty, -c]}(A_y) = 1-q_y-q^+_y
\]

are norm continuous on $U'$.
Let $A^-_y$ and $A^+_y$ be the restrictions of $A_y$ to the ranges of
projections $q^-_y$ and $q^+_y$ respectively.
Then 
\begin{equation}\label{eq:AAA}
  \f(A_y) = \f(A^-_y)+\f(A'_y)+\f(A^+_y). 
\end{equation}
Since $\f(c)\in(1-\delta,1)$, we have the inequalities 
\begin{equation}\label{eq:fA}
	\norm{\f(A^+_y) - q^+_y} < \delta \quad \text{and} \quad \norm{\f(A^-_y) + q^-_y} < \delta  
\end{equation}
for every $y\in U'$. 
Let $U''\subset U'$ be a neighborhood of $x$ such that 
\begin{equation}\label{eq:fA'}
 \norm{\f(A'_y)-\f(A'_x)}<\delta, \quad \norm{q_y^- - q_x^-}<\delta, \;\text{and } \norm{q_y^+ - q_x^+}<\delta 
\end{equation}
for every $y\in U''$.
Combining \eqref{eq:AAA}, \eqref{eq:fA}, and \eqref{eq:fA'}, we get 
\[
\norm{\f(A_y)-\f(A_x)} < 7\delta
\] 
for every $y\in U''$. 
Since $x\in U$ and $\delta>0$ were chosen arbitrarily, 
$\f(\A)$ is norm continuous and thus $\A$ is Riesz continuous in a given local trivialization.
\endproof

\sub{Compactly-polarized operators.}\hspace*{-0.2em}
A \emph{compactly-polarized operator} in a Hilbert space $H$ is defined in \cite{i2}
as an essentially unitary self-adjoint operator in $H$ with the norm $\leq 1$
and the essential spectrum $\set{-1, 1}$ 
or, what is the same, as an operator from Atiyah--Singer \cite{as} space $\widehat{F}_*$ related to $H$. 

Let $\A$ be a family of compactly-polarized operators in $\H$.
Such a family is never a discrete-spectrum family because neither of the families $\A-1$, $\A+1$ is Fredholm. 
We will say that it is a weakly discrete-spectrum family if $\A-\lambda$ is a Fredholm family for every $\lambda\in(-1,1)$.
Clearly, a family of self-adjoint operators is discrete-spectrum if and only if 
its bounded transform is weakly discrete-spectrum.

The following two statements are proved in exactly the same manner as Lemma \ref{lem:1} and Theorem \ref{thm:2},
so we omit the proofs.
Theorem \ref{thm:3} strengthens Lemma 5.1 from \cite{i2}.

\begin{lem}\label{lem:2}
%Let $\A$ be a family of self-adjoint operators with the norm $\leq 1$.
Let $\A$ be a family of compactly-polarized operators.
Then $\A$ is a weakly discrete-spectrum family if and only if
for every $x\in X$ and $0<b<1$ there is a pair $(U,c)$ adapted to $\A$ such that $x\in U$ and $c>b$.
\end{lem}

\begin{thm}\label{thm:3}
Let $\A$ be a weakly discrete-spectrum family.
Then every local trivialization strictly adapted to $\A$ is fully adapted to $\A$.
In particular, every weakly discrete-spectrum and strictly Fredholm family is fully Fredholm.
\end{thm}

%\begin{small}

%\end{small}

\end{document}